
\input amssym.def
\input amssym.tex

\hsize=13.50cm    
\vsize=18cm       
\parindent=12pt   \parskip=0pt      
\pageno=1 


\hoffset=15mm    
\voffset=1cm    
 

\ifnum\mag=\magstep1
\hoffset=-2mm   
\voffset=.8cm   
\fi


\pretolerance=500 \tolerance=1000  \brokenpenalty=5000

\catcode`\@=11

\font\eightrm=cmr8         \font\eighti=cmmi8
\font\eightsy=cmsy8        \font\eightbf=cmbx8
\font\eighttt=cmtt8        \font\eightit=cmti8
\font\eightsl=cmsl8        \font\sixrm=cmr6
\font\sixi=cmmi6           \font\sixsy=cmsy6
\font\sixbf=cmbx6


\font\tengoth=eufm10       \font\tenbboard=msbm10
\font\eightgoth=eufm10 at 8pt      \font\eightbboard=msbm10 at 8 pt
\font\sevengoth=eufm7      \font\sevenbboard=msbm7
\font\sixgoth=eufm7 at 6 pt        \font\fivegoth=eufm5

 \font\tencyr=wncyr10       
\font\eightcyr=wncyr10 at 8 pt      
\font\sevencyr=wncyr10 at 7 pt      
\font\sixcyr=wncyr10 at 6 pt


\skewchar\eighti='177 \skewchar\sixi='177
\skewchar\eightsy='60 \skewchar\sixsy='60


\newfam\gothfam           \newfam\bboardfam
\newfam\cyrfam

\def\tenpoint{%
  \textfont0=\tenrm \scriptfont0=\sevenrm \scriptscriptfont0=\fiverm
  \def\rm{\fam\z@\tenrm}%
  \textfont1=\teni  \scriptfont1=\seveni  \scriptscriptfont1=\fivei
  \def\oldstyle{\fam\@ne\teni}\let\old=\oldstyle
  \textfont2=\tensy \scriptfont2=\sevensy \scriptscriptfont2=\fivesy
  \textfont\gothfam=\tengoth \scriptfont\gothfam=\sevengoth
  \scriptscriptfont\gothfam=\fivegoth
  \def\goth{\fam\gothfam\tengoth}%
  \textfont\bboardfam=\tenbboard \scriptfont\bboardfam=\sevenbboard
  \scriptscriptfont\bboardfam=\sevenbboard
  \def\bb{\fam\bboardfam\tenbboard}%
 \textfont\cyrfam=\tencyr \scriptfont\cyrfam=\sevencyr
  \scriptscriptfont\cyrfam=\sixcyr
  \def\cyr{\fam\cyrfam\tencyr}%
  \textfont\itfam=\tenit
  \def\it{\fam\itfam\tenit}%
  \textfont\slfam=\tensl
  \def\sl{\fam\slfam\tensl}%
  \textfont\bffam=\tenbf \scriptfont\bffam=\sevenbf
  \scriptscriptfont\bffam=\fivebf
  \def\bf{\fam\bffam\tenbf}%
  \textfont\ttfam=\tentt
  \def\tt{\fam\ttfam\tentt}%
  \abovedisplayskip=12pt plus 3pt minus 9pt
  \belowdisplayskip=\abovedisplayskip
  \abovedisplayshortskip=0pt plus 3pt
  \belowdisplayshortskip=4pt plus 3pt 
  \smallskipamount=3pt plus 1pt minus 1pt
  \medskipamount=6pt plus 2pt minus 2pt
  \bigskipamount=12pt plus 4pt minus 4pt
  \normalbaselineskip=12pt
  \setbox\strutbox=\hbox{\vrule height8.5pt depth3.5pt width0pt}%
  \let\bigf@nt=\tenrm       \let\smallf@nt=\sevenrm
  \normalbaselines\rm}

\def\eightpoint{%
  \textfont0=\eightrm \scriptfont0=\sixrm \scriptscriptfont0=\fiverm
  \def\rm{\fam\z@\eightrm}%
  \textfont1=\eighti  \scriptfont1=\sixi  \scriptscriptfont1=\fivei
  \def\oldstyle{\fam\@ne\eighti}\let\old=\oldstyle
  \textfont2=\eightsy \scriptfont2=\sixsy \scriptscriptfont2=\fivesy
  \textfont\gothfam=\eightgoth \scriptfont\gothfam=\sixgoth
  \scriptscriptfont\gothfam=\fivegoth
  \def\goth{\fam\gothfam\eightgoth}%
  \textfont\cyrfam=\eightcyr \scriptfont\cyrfam=\sixcyr
  \scriptscriptfont\cyrfam=\sixcyr
  \def\cyr{\fam\cyrfam\eightcyr}%
  \textfont\bboardfam=\eightbboard \scriptfont\bboardfam=\sevenbboard
  \scriptscriptfont\bboardfam=\sevenbboard
  \def\bb{\fam\bboardfam}%
  \textfont\itfam=\eightit
  \def\it{\fam\itfam\eightit}%
  \textfont\slfam=\eightsl
  \def\sl{\fam\slfam\eightsl}%
  \textfont\bffam=\eightbf \scriptfont\bffam=\sixbf
  \scriptscriptfont\bffam=\fivebf
  \def\bf{\fam\bffam\eightbf}%
  \textfont\ttfam=\eighttt
  \def\tt{\fam\ttfam\eighttt}%
  \abovedisplayskip=9pt plus 3pt minus 9pt
  \belowdisplayskip=\abovedisplayskip
  \abovedisplayshortskip=0pt plus 3pt
  \belowdisplayshortskip=3pt plus 3pt 
  \smallskipamount=2pt plus 1pt minus 1pt
  \medskipamount=4pt plus 2pt minus 1pt
  \bigskipamount=9pt plus 3pt minus 3pt
  \normalbaselineskip=9pt
  \setbox\strutbox=\hbox{\vrule height7pt depth2pt width0pt}%
  \let\bigf@nt=\eightrm     \let\smallf@nt=\sixrm
  \normalbaselines\rm}

\tenpoint


\def\pc#1{\bigf@nt#1\smallf@nt}         \def\pd#1 {{\pc#1} }


\catcode`\;=\active
\def;{\relax\ifhmode\ifdim\lastskip>\z@\unskip\fi
\kern\fontdimen2  -1.2 \fontdimen3 \string;}

\catcode`\:=\active
\def:{\relax\ifhmode\ifdim\lastskip>\z@\unskip\fi\penalty\@M\ \fi\string:}

\catcode`\!=\active
\def!{\relax\ifhmode\ifdim\lastskip>\z@
\unskip\fi\kern\fontdimen2  -1.1 \fontdimen3 \string!}

\catcode`\?=\active
\def?{\relax\ifhmode\ifdim\lastskip>\z@
\unskip\fi\kern\fontdimen2  -1.1 \fontdimen3 \string?}

\def\^#1{\if#1i{\accent"5E\i}\else{\accent"5E #1}\fi}
\def\"#1{\if#1i{\accent"7F\i}\else{\accent"7F #1}\fi}

\frenchspacing


\newtoks\auteurcourant      \auteurcourant={\hfil}
\newtoks\titrecourant       \titrecourant={\hfil}

\newtoks\hautpagetitre      \hautpagetitre={\hfil}
\newtoks\baspagetitre       \baspagetitre={\hfil}

\newtoks\hautpagegauche     
\hautpagegauche={\eightpoint\rlap{\folio}\hfil\the\auteurcourant\hfil}
\newtoks\hautpagedroite     
\hautpagedroite={\eightpoint\hfil\the\titrecourant\hfil\llap{\folio}}

\newtoks\baspagegauche      \baspagegauche={\hfil} 
\newtoks\baspagedroite      \baspagedroite={\hfil}

\newif\ifpagetitre          \pagetitretrue  


\headline={\ifpagetitre\the\hautpagetitre
\else\ifodd\pageno\the\hautpagedroite\else\the\hautpagegauche\fi\fi}

\footline={\ifpagetitre\the\baspagetitre\else
\ifodd\pageno\the\baspagedroite\else\the\baspagegauche\fi\fi
\global\pagetitrefalse}


\def\raggedbottom{\topskip 10pt plus 36pt\r@ggedbottomtrue}



\def\pointir{\unskip . --- \ignorespaces}


\def\Bigbreak{\vskip-\lastskip\bigbreak}
\def\Medbreak{\vskip-\lastskip\medbreak}


\def\ctexte#1\endctexte{%
  \hbox{$\vcenter{\halign{\hfill##\hfill\crcr#1\crcr}}$}}


\long\def\ctitre#1\endctitre{%
    \ifdim\lastskip<24pt\vskip-\lastskip\bigbreak\bigbreak\fi
  		\vbox{\parindent=0pt\leftskip=0pt plus 1fill
          \rightskip=\leftskip
          \parfillskip=0pt\bf#1\par}
    \bigskip\nobreak}

\long\def\section#1\endsection{%
\vskip 0pt plus 3\normalbaselineskip
\penalty-250
\vskip 0pt plus -3\normalbaselineskip
\Bigbreak
\message{[section \string: #1]}{\bf#1\unskip}\pointir}

\long\def\sectiona#1\endsection{%
\vskip 0pt plus 3\normalbaselineskip
\penalty-250
\vskip 0pt plus -3\normalbaselineskip
\Bigbreak
\message{[sectiona \string: #1]}%
{\bf#1}\medskip\nobreak}

\long\def\subsection#1\endsubsection{%
\Medbreak
{\it#1\unskip}\pointir}

\long\def\subsectiona#1\endsubsection{%
\Medbreak
{\it#1}\par\nobreak}

\def\rem#1\endrem{%
\Medbreak
{\it#1\unskip} : }

\def\remp#1\endrem{%
\Medbreak
{\pc #1\unskip}\pointir}

\def\rema#1\endrem{%
\Medbreak
{\it #1}\par\nobreak}

\def\newparwithcolon#1\endnewparwithcolon{
\Medbreak
{#1\unskip} : }

\def\newparwithpointir#1\endnewparwithpointir{
\Medbreak
{#1\unskip}\pointir}

\def\newpara#1\endnewpar{
\Medbreak
{#1\unskip}\smallskip\nobreak}


\long\def\th#1 #2\enonce#3\endth{%
   \Medbreak
   {\pc#1} {#2\unskip}\pointir{\it #3}\medskip}

\long\def\tha#1 #2\enonce#3\endth{%
   \Medbreak
   {\pc#1} {#2\unskip}\par\nobreak{\it #3}\medskip}


\long\def\Th#1 #2 #3\enonce#4\endth{%
   \Medbreak
   #1 {\pc#2} {#3\unskip}\pointir{\it #4}\medskip}

\long\def\Tha#1 #2 #3\enonce#4\endth{%
   \Medbreak
   #1 {\pc#2} #3\par\nobreak{\it #4}\medskip}


\def\decale#1{\smallbreak\hskip 28pt\llap{#1}\kern 5pt}
\def\decaledecale#1{\smallbreak\hskip 34pt\llap{#1}\kern 5pt}
\def\puce{\smallbreak\hskip 6pt{$\scriptstyle\bullet$}\kern 5pt}



\def\displaylinesno#1{\displ@y\halign{
\hbox to\displaywidth{$\@lign\hfil\displaystyle##\hfil$}&
\llap{$##$}\crcr#1\crcr}}


\def\ldisplaylinesno#1{\displ@y\halign{ 
\hbox to\displaywidth{$\@lign\hfil\displaystyle##\hfil$}&
\kern-\displaywidth\rlap{$##$}\tabskip\displaywidth\crcr#1\crcr}}


\def\eqalign#1{\null\,\vcenter{\openup\jot\m@th\ialign{
\strut\hfil$\displaystyle{##}$&$\displaystyle{{}##}$\hfil
&&\quad\strut\hfil$\displaystyle{##}$&$\displaystyle{{}##}$\hfil
\crcr#1\crcr}}\,}


\def\system#1{\left\{\null\,\vcenter{\openup1\jot\m@th
\ialign{\strut$##$&\hfil$##$&$##$\hfil&&
        \enskip$##$\enskip&\hfil$##$&$##$\hfil\crcr#1\crcr}}\right.}


\let\@ldmessage=\message

\def\message#1{{\def\pc{\string\pc\space}%
                \def\'{\string'}\def\`{\string`}%
                \def\^{\string^}\def\"{\string"}%
                \@ldmessage{#1}}}



\def\up#1{\raise 1ex\hbox{\smallf@nt#1}}


\def\qed{\raise -2pt\hbox{\vrule\vbox to 10pt{\hrule width 4pt
                 \vfill\hrule}\vrule}}

\def\cqfd{\unskip\penalty 500\quad\qed\medbreak}

\def\virg{\raise .4ex\hbox{,}}   


\def\build#1_#2^#3{\mathrel{
\mathop{\kern 0pt#1}\limits_{#2}^{#3}}}


\def\boxit#1#2{%
\setbox1=\hbox{\kern#1{#2}\kern#1}%
\dimen1=\ht1 \advance\dimen1 by #1 \dimen2=\dp1 \advance\dimen2 by #1 
\setbox1=\hbox{\vrule height\dimen1 depth\dimen2\box1\vrule}%
\setbox1=\vbox{\hrule\box1\hrule}%
\advance\dimen1 by .6pt \ht1=\dimen1 
\advance\dimen2 by .6pt \dp1=\dimen2  \box1\relax}


\catcode`\@=12

\showboxbreadth=-1  \showboxdepth=-1



\magnification=\magstep1
\hsize=17,5truecm  
\vsize=25.5truecm 
\hoffset=-0.9truecm 
\voffset=-0.8truecm
\topskip=1truecm
\footline={\tenrm\hfil\folio\hfil}
\raggedbottom
\abovedisplayskip=3mm 
\belowdisplayskip=3mm
\abovedisplayshortskip=0mm
\belowdisplayshortskip=2mm
\normalbaselineskip=12pt  
\normalbaselines

\def\deg{{\rm deg}}

\def\G{{\bf G}}
\def\A{{\bf A}}
\def\I{{\bf I}}
\def\F{{\bf F}}

\def\Z{{\bf Z}}

\def\P{{\bf P}}
\def\hom{{\rm Hom}}
\def\deg{{\rm deg}}
\

\bigskip

{\bf  Quelques questions d'approximation faible pour les tores alg\'ebriques}

\bigskip

J.-L. Colliot-Th\'el\`ene et V. Suresh

\bigskip

{\bf Introduction}

\medskip
Il est bien connu que les tores alg\'ebriques d\'efinis sur un corps global ne satisfont
pas n\'ecessairement l'approximation faible : \'etant donn\'e un ensemble $S$ fini non vide de places
du corps global $K$, et un $K$-tore $T$, le groupe des points rationnels $T(K)$
n'est pas n\'ecessairement dense dans le produit $\prod_{v \in S} T(K_{v})$.

Pour chaque place $v$, le groupe  $T(K_{v})$ est ici muni de la topologie induite par celle du corps local $K_{v}$, compl\'et\'e de $K$ en la place $v$. C'est un groupe topologique commutatif localement compact. Notons $T(O_{v}) \subset T(K_{v})$ son sous-groupe compact maximal.  Dans plusieurs contextes, on a \'et\'e amen\'e \`a se poser la question d'approximation suivante,  o\`u l'on demande moins que l'approximation faible.

\medskip

{\bf Question semi-locale } {\it Le sous-groupe ouvert $T(K).\prod_{v \in S} T(O_{v})$ de $\prod_{v \in S} T(K_{v})$
co\"{\i}ncide-t-il avec $\prod_{v \in S} T(K_{v})$ ?
}

\medskip

En d'autres termes, l'application naturelle de $T(K) $ vers le groupe discret $\prod_{v \in S} T(K_{v})/ T(O_{v}) $ est-elle surjective ? Lorsque l'ensemble $S$ est r\'eduit \`a une place, la question fut pos\'ee par Bruhat et Tits (voir [CTS2], Remark 8.3 p. 192). 

Dans cet article, sur un corps global de caract\'eristique positive, nous r\'epondons n\'egativement \`a la question semi-locale  (mais laissons ouverte
la question de Bruhat et Tits). Nous r\'epondons aussi n\'egativement \`a la question purement
locale suivante ([CTS2],  Remark 8.3 p. 192). 

\medskip

{\bf Question locale } {\it  Soient $K$ un corps local, $T$ un $K$-tore, $T(O_{K}) \subset T(K)$ le sous-groupe compact
maximal, $RT(K) \subset T(K)$ le sous-groupe des \'el\'ements R-\'equivalents \`a l'\'el\'ement neutre.
A-t-on $T(O_{K}).RT(K)=T(K)$ ?}

\medskip

Enfin, lorsque $K$ est un corps de fonctions d'une variable sur un corps fini,
nous r\'epondons n\'egativement \`a une question soulev\'ee par D. Bourqui.
Cette question ({\bf Question globale}, \S 3 ci-dessous) est apparue naturellement dans le travail
[B] sur la fonction z\^{e}ta des hauteurs sur une vari\'et\'e torique sur un corps $K$
de fonctions d'une variable sur un corps fini.
La r\'eponse n\'egative que nous apportons permet \`a Bourqui de montrer 
que la constante d\'efinie par Peyre [P1]  et  Batyrev/Tschinkel [BT] (voir le rapport [P2])
pour les vari\'et\'es toriques sur   un corps de nombres $K$
doit, dans le cas fonctionnel, \^etre multipli\'ee par une certaine constante (\`a valeurs enti\`eres)
non n\'ecessairement \'egale \`a $1$.

On trouvera au  \S 1 des rappels de [CTS1]. Un bref \S 2 discute les questions locale et semi-locale.
Au \S 3, on pr\'esente la question globale.  Le \S 4 contient la description alg\'ebrique
du tore que nous utilisons pour donner des contre-exemples. Le \S 5 contient
la r\'eponse n\'egative \`a la question locale. Le \S 6 contient la r\'eponse n\'egative
aux questions semi-locale et globale.

\bigskip

{\bf \S 1 R\'esolutions flasques et coflasques, R-\'equivalence : rappels}

\medskip

Soit  $L/K$ une extension finie de corps,
galoisienne de groupe $G$.
Etant donn\'e un $K$-tore $T$ d\'eploy\'e par $L$,
c'est-\`a-dire un $K$-groupe alg\'ebrique $T$ tel que le $L$-groupe
$T\times_{K}L$ est $L$-isomorphe \`a un produit de groupes multiplicatifs $\G_{m,L}$,
on note $T^*={\rm Hom}_{L-{\rm groupe}}(T_{L}, \G_{m,L})$ son groupe des caract\`eres.
C'est un $G$-r\'eseau. On note $T_{*}={\rm Hom}_{L-{\rm groupe}}(\G_{m,L}, T_{L})$
le groupe des cocaract\`eres de $T\times_KL$. C'est le $G$-r\'eseau dual du $G$-r\'eseau $T^*$,
c'est-\`a-dire que l'on a $T^*={\rm Hom}_{\Z}(T_{*},\Z)$.

On sait (Endo-Miyata, Voskresenski\u{\i}, voir [CTS1] \S 1, Lemme 3) 
que pour tout
$G$-r\'eseau
$T_*$ on peut trouver une suite exacte de $G$-r\'eseaux
$$0 \to F_* \to P_* \to T_* \to 0 \hskip3cm (1)$$
avec $P_*$ un $G$-module de permutation
et $F_*$ un $G$-module coflasque, c'est-\`a-dire
tel que $H^1(H,F_*)=0$ pour tout sous-groupe $H \subset G$.
Une telle suite est dite {\it r\'esolution coflasque} du $G$-r\'eseau $T_*$.
Si 
$$0 \to F_{1*} \to P_{1*} \to T_* \to 0,$$
et
$$0 \to F_{2*} \to P_{2*} \to T_* \to 0,$$
sont deux telles r\'esolutions, on montre ([CTS1], \S 1, Lemme 5; [ CTS2], Lemma 0.6) qu'il existe 
un isomorphisme de $G$-r\'eseaux $
F_{1*} \oplus P_{2*} \simeq F_{2*} \oplus P_{1*}$.
Plus pr\'ecis\'ement, si l'on note $M_{*}$ le $G$-r\'eseau produit fibr\'e
de $ P_{1*} \to T_* $ et $P_{2*} \to T_*$, les projections $M_{*} \to P_{1*}$
et $M_{*} \to P_{2*}$ sont $G$-scind\'ees ([CTS1], Lemmes 1 et 5).

La suite (1) induit une suite exacte de $K$-tores d\'eploy\'es par $L$
$$ 1 \to F \to P \to T \to 1 \hskip3cm (2)$$
avec $P$ un $K$-tore quasitrivial et $F$ un $K$-tore flasque.
Une telle suite est appel\'ee une {\it r\'esolution flasque} du $K$-tore $T$.

Rappelons ici que le $G$-module $T(L)$ des $L$-points d'un $K$-tore $T$
d\'eploy\'e par $L$ est le $G$-module $T_* \otimes L^{\times}=T_*
\otimes_{\Z} L^{\times}$
\'equip\'e de l'action diagonale de $G$.

Deux points $p,q$ de $T(K)$ sont dits R-\'equivalents s'il existe
un ouvert $U$ de la droite projective $\P^1_{k}$ et
un $K$-morphisme  $\varphi : U \to T$ tels que $p,q \in \varphi(U(K))$
(on d\'emontre que c'est une relation d'\'equivalence).
Comme il est \'etabli au \S 5 de [CTS1], la suite exacte
$$ P(K) \to T(K) \to H^1(K,F) \to 0$$
tir\'ee de (2) par cohomologie galoisienne calcule la $R$-\'equivalence
sur le groupe $T(K)=(T_* \otimes L^{\times})^G$ des points $K$-rationnels
du tore $T$.  L'image de $P(K)$ dans $T(K)$ est
exactement le sous-groupe $RT(K) \subset T(K)$ des points
$R$-\'equivalents \`a l'\'el\'ement neutre dans $T(K)$. En d'autres termes,
$T(K)/R = H^1(K,F)$.

On sait (Endo-Miyata, cf. [CTS1], Prop. 2 p. 184) que lorsque le groupe $G$ 
est m\'etacyclique, i.e. a tous ses sous-groupes
de Sylow cycliques, tout $G$-module coflasque $F_{*}$ est facteur direct d'un $G$-module de permutation.
Ceci implique $H^1(K,F)=0$. Ainsi, si un $K$-tore $T$ est d\'eploy\'e par
une extension $L/K$ m\'etacyclique, on a $T(K)/R=1$.

\bigskip

{\bf \S 2 La question locale et la question semi-locale}

\medskip

Soient $K$ un corps local non archim\'edien et $L/K$ une extension finie galoisienne
de groupe de Galois $G$. Soit $O_{K}$, resp. $O_{L}$, 
l'anneau des entiers de $K$, resp. $L$. 
La valuation normalis\'ee $v_{L} : L^* \to \Z$ donne naissance \`a la suite exacte de $G$-modules
$$ 1 \to O_{L}^{\times} \to L^{\times} \to \Z \to 0,$$
o\`u l'action de $G$ sur $\Z$ est triviale, la fl\`eche $L^{\times} \to \Z $ \'etant donn\'ee par la valuation
(normalis\'ee) $v_{L}$ de $L$. 

Soit $T$ un $K$-tore d\'eploy\'e par $L$, $T_{*}$ son groupe 
des cocaract\`eres. C'est un $G$-r\'eseau.

De la suite exacte de la valuation normalis\'ee sur $L$  on d\'eduit la suite exacte de $G$-modules
$$ 0 \to T_{*} \otimes O_{L}^{\times} \to T_{*} \otimes L^{\times} \to T_{*} \to 0.$$

Le groupe $(T_{*} \otimes O_{L}^{\times})^G \subset (T_{*} \otimes {L^{\times}})^G=T(K)$
est {\it not\'e}  $T(O_{K})$. C'est le sous-groupe compact maximal de $T(K)$.

Le $K$-tore $T$ est dit anisotrope si l'une des conditions suivantes est satisfaite :
$T_{*} ^G=0$, ou ${T^{*}}^G=0$. Si $T$ est anisotrope, alors $T(O_{K})=T(K)$ et  
$T(K)$ est compact  (comme il est bien connu, et comme il est facile \`a \'etablir \`a partir des
suites ci-dessus, cette condition n\'ecessaire d'anisotropie est une condition suffisante.)

\bigskip

Commen\c cons par commenter la question locale.

\medskip

{\bf Question locale} {\it  Soit $K$ un corps local, $T$ un $K$-tore, $T(O_{K}) \subset T(K)$ le sous-groupe compact
maximal, $RT(K) \subset T(K)$ le sous-groupe des \'el\'ements R-\'equivalents \`a l'\'el\'ement neutre.
A-t-on $T(O_{K}).RT(K)=T(K)$ ?}

En d'autres termes, tout
\'el\'ement de $T(K)$ est-il produit d'un \'el\'ement de $RT(K)$ et d'un
\'el\'ement de $T(O_K)$ ?
En d'autres termes encore, le sous-groupe compact maximal $T(O_K)$
rencontre-t-il toutes les classes pour la $R$-\'equivalence sur $T(K)$ ?

La r\'eponse est trivialement positive si  $T(K)/R=1$.  Elle est   positive
dans de nombreux cas.

\bigskip

{\bf Proposition 2.1} 
{\it  La question locale a une r\'eponse affirmative dans chacun des cas suivants.

(i) Le $K$-tore $T$ a bonne r\'eduction.

(ii) Le $K$-tore $T$ est d\'eploy\'e par une extension m\'etacyclique.

(iii) Le $K$-tore $T$ est anisotrope.

(iv) Le $K$-tore $T$ est d\'eploy\'e par une extension finie galoisienne $L/K$ et
admet une r\'esolution flasque du type
$$1 \to F \to  (R_{L/K}\G_{m})^n \to T \to 1.$$

(v) (Bourqui) Le $K$-tore $T$ est d\'eploy\'e par une extension $L/K$ totalement ramifi\'ee.}

\medskip

{\it D\'emonstration} 

Dans le cas (i), le tore est d\'eploy\'e par une extension cyclique, ce cas est
un cas particulier de (ii).  Dans le cas (ii), on a $T(K)/R=1$ comme il a \'et\'e rappel\'e
au \S 1. Ces cas sont donc \'evidents. Le cas (iii) l'est aussi, car, comme il a \'et\'e rappel\'e ci-dessus, si $T$ est anisotrope, alors  $T(O_K)=T(K)$.

Etablissons le cas (iv).  Soit $G$ le groupe de Galois de $L/K$.
Le $K$-tore $R_{L/K}\G_{m}$ a le module galoisien $\Z[G]$ pour groupe des  cocaract\`eres. 
Soit 
$$0 \to F_* \to P_* \to T_* \to 0 $$
une suite exacte de $G$-r\'eseaux du type (1), avec $P_{*} = (\Z[G])^n$ pour $n>0$ convenable.

En tensorisant la suite de type (1) (qui est $\Z$-scind\'ee) par $O_{L}^{\times}$, et en prenant la $G$-cohomologie de la suite exacte obtenue, on obtient la suite exacte
$$ (T_{*} \otimes O_{L}^{\times})^G  \to H^1(G,F_{*}\otimes O_{L}^{\times}) \to H^1(G,P_{*}\otimes O_{L}^{\times}),$$
soit encore
$$T(O_{K}) \to H^1(G,F_{*}\otimes O_{L}^{\times}) \to 0,$$ 
tout groupe de la forme $H^r(G,\Z[G]\otimes  M)$ avec $r>0$ \'etant nul.
Si l'on tensorise la suite de la valuation par le groupe ab\'elien libre  $F_*$, et si l'on prend
la suite de cohomologie de la suite exacte courte de $G$-modules ainsi obtenue, on obtient la suite
exacte
$$ H^1(G, F_{*}\otimes O_{L}^{\times}) \to H^1(G, F_{*}\otimes L^{\times})
\to H^1(G, F_{*}).$$
Le dernier groupe est nul, car $F_{*}$ est coflasque. Ainsi la fl\`eche compos\'ee
$$T(O_{K}) \to H^1(G,F_{*}\otimes O_{L}^{\times}) \to H^1(G, F_{*}\otimes L^{\times})$$
est surjective. Comme cette fl\`eche co\"{\i}ncide avec la fl\`eche compos\'ee
$T(O_{K}) \to T(K) \to T(K)/R$, ceci \'etablit l'assertion dans le cas (iv).

On notera que les $K$-tores normiques $R^1_{L/K}\G_{m} = {\rm Ker} [N_{L/K} : R_{L/K}\G_{m} \to \G_{m}]$, pour $L/K$ extension finie galoisienne, sont du type (iv) ([CTS1], \S 6, Prop.~15 p.~206).

Consid\'erons maintenant le cas (v), qui nous a  \'et\'e signal\'e par D. Bourqui.
Le groupe $T(O_{K})$ est le noyau de la fl\`eche $T(K)=(T_*\otimes L^{\times})^G \to  {T_*}^G$
induite par la valuation (normalis\'ee) $v_{L } : L^{\times} \to \Z$. L'accouplement naturel non d\'eg\'en\'er\'e
$ T_{*}  \times  T^* \to \Z$ induit 
un homomorphisme  ${T_*}^G \to {\rm Hom}({T^*}^G, \Z)$ dont on v\'erifie qu'il est injectif
(Lemme 3.1 ci-apr\`es). Ainsi $T(O_{K})$ est le noyau de l'application compos\'ee 
$$\phi_{K,L } : T(K)=(T_*\otimes L^{\times})^G \to  {T_*}^G \to {\rm Hom}({T^*}^G, \Z),$$
o\`u la premi\`ere f\l\`eche est induite par la valuation $v_{L}$.

Un \'el\'ement de ${T^*}^G$  correspond \`a un $K$-homomorphisme $T \to \G_{m,K}$.
Un tel $K$-morphisme induit un homomorphisme $T(K) \to K^{\times}$ que l'on peut composer
avec la valuation (normalis\'ee) $v_{K }: K^{\times} \to \Z$. Ceci d\'efinit un homomorphisme
$$\psi_{K} : T(K) \to {\rm Hom}({T^*}^G, \Z).$$
Soit $e$ l'indice de ramification de  $L$ sur  $K$. On v\'erifie ais\'ement la formule 
$\phi_{K,L}=e \psi_{K},$ qui implique en particulier  que le noyau de $\psi_{K}$ est $T(O_{K})$ (le groupe  ${\rm Hom}({T^*}^G, \Z)$ est sans torsion).
L'application $T(L) \to {\rm Hom}({T^*},\Z)$ induite par la valuation sur $L$ est clairement surjective.
L'application de restriction ${\rm Hom}({T^*},\Z) \to {\rm Hom}({T^*}^G,\Z)$ est surjective, car
le groupe ab\'elien ${T^*}^G$ est facteur direct dans ${T^*}$ (le quotient \'etant sans torsion).
Ainsi l'application compos\'ee
 $$T(L) \to {\rm Hom}({T^*},\Z) \to  {\rm Hom}({T^*}^G,\Z)$$
est surjective.  La compos\'ee de cette application avec l'inclusion $T(K) \subset T(L)$
est $\phi_{K,L}$.
Cette application est $G$-\'equivariante, l'action de $G$ sur 
le groupe ${\rm Hom}({T^*}^G,\Z)$ \'etant l'action triviale.
Soit $\theta \in {\rm Hom}({T^*}^G,\Z)$. Soit $\beta \in T(L)$ d'image $\theta$
par l'application ci-dessus. L'image de $\alpha=\prod_{g \in G} g.\beta$ est 
$[L:K]\theta$.  On a donc
$$e\psi_{K}(\alpha)= \phi_{K,L}(\alpha)= [L:K] \theta \in {\rm Hom}({T^*}^G,\Z).$$

Supposons l'extension $L/K$  totalement ramifi\'ee, i.e. $e=[L:K]$. Alors
$\psi_{K}(\alpha)=\theta $ dans le groupe  ab\'elien libre ${\rm Hom}({T^*}^G,\Z).$ Comme $\alpha$
est la norme de $\beta \in T(L)$, ceci \'etablit $$\psi_{K}(N_{L/K}(T(L)))={\rm Hom}({T^*}^G,\Z).$$ 
Soit $1 \to F \to P \to T \to 1$
une r\'esolution flasque du $K$-tore $T$ par des $K$-tores d\'eploy\'es par $L$.
L'homomorphisme induit $P(L) \to T(L)$ est surjectif, l'application compos\'ee
$P(L) \to T(L) \to T(K) \to {\rm Hom}({T^*}^G,\Z)$ l'est donc aussi, o\`u $T(L) \to T(K)$
est la norme. Ceci implique que l'application compos\'ee $P(K) \to T(K) \to 
{\rm Hom}({T^*}^G,\Z)$ est surjective.  La fl\`eche $T(K) \to 
{\rm Hom}({T^*}^G,\Z)$ est ici $\psi_{K}$, son noyau est $T(O_{K})$. On voit donc
que $T(K)$ est engendr\'e par $T(O_{K})$ et l'image de $P(K) \to T(K)$,
qui est le sous-groupe $RT(K)$.
\cqfd

\bigskip

Discutons maintenant la {\bf question semi-locale}.
Soient  
 $K$ un corps global et $S$ un ensemble fini non vide de places de $K$.
Soit $T$ un $K$-tore d\'eploy\'e par une extension finie galoisienne $L/K$ de groupe 
de Galois $G$. Soit 
$$0 \to F_* \to P_* \to T_* \to 0 $$
une suite exacte de $G$-r\'eseaux du type (1),
induisant une suite exacte de $K$-tores
$$1\to F  \to P  \to T  \to 1.$$

On a les inclusions \'evidentes suivantes :

$$ T(K).\prod_{v \in S} T(O_{v}) \subset T(K).\prod_{v \in S} (T(O_{v}).RT(K_{v}))
\subset \prod_{v \in S} T(K_{v}).$$

Le premier groupe
est un sous-groupe ouvert de $ \prod_{v \in S} T(K_{v})$ contenant $T(K)$.
Le $K$-tore $P$ est quasi-trivial, donc est un ouvert de Zariski d'un espace affine.  Ainsi
$P(K)$ est dense dans $\prod_{v \in S} P(K_{v})$. 
Ceci implique que l'image de $P(K)$ dans $T(K)$  est dense
dans le produit des images des $P(K_{v})$ dans  $\prod_{v \in S} T(K_{v})$,
c'est-\`a-dire dans $\prod_{v \in S}  RT(K_{v})$. Tout point de
$\prod_{v \in S}  RT(K_{v})$ peut donc s'\'ecrire comme le produit
d'un \'el\'ement d'un \'el\'ement de $T(K)$ (dans l'image de $P(K)$) et
d'un \'el\'ement de l'ouvert $\prod_{v \in S} T(O_{v})$.
Ainsi {\it la premi\`ere inclusion ci-dessus est une \'egalit\'e}.

\medskip

Ceci permet de reformuler la {\bf question semi-locale} de la fa\c con suivante :
\medskip

{\it L'application naturelle 
$ T(K).\prod_{v \in S} T(O_{v}) \to \prod_{v \in S} H^1(K_v,F)$
est-elle surjective ?}

\medskip

Ceci montre aussi : 

\medskip

{\bf Proposition 2.2} 
{\it Une r\'eponse affirmative \`a la question locale (pour chaque $K_{v}$-tore
$T\times_KK_{v}$)  implique  une r\'eponse affirmative \`a la question semi-locale.}\cqfd

\bigskip

{\bf \S 3 La question globale (cas fonctionnel)}

\medskip

Soient $\F$ un corps fini et $K$ un corps de fonctions d'une variable
sur le corps $\F$, c'est-\`a-dire une extension de type fini, de
degr\'e de transcendance un, du corps $\F$. On ne suppose pas le corps
$\F$ alg\'ebriquement ferm\'e dans $K$. 
Soit $L/K$ une extension finie de corps. Soit $\Omega_L$ 
l'ensemble des places de $L$, et pour $w \in
\Omega_L$, soient $L_w$ le compl\'et\'e de $L$ en $w$ et $O_w$ son anneau
des entiers. 
On note encore $w : L_w^{\times} \to \Z$ la valuation normalis\'ee
(i.e. d'image le groupe $\Z$ tout entier).
Le corps
r\'esiduel  $\F_w$ du corps local $L_w$ est une extension finie de $\F$.
Soit 
$\I_L$ le groupe des id\`eles de $L$, c'est-\`a-dire le
produit restreint des $L_{w}^{\times}$ pour $w \in \Omega_L$.
On note $$\deg_{L,\F} : \I_L \to \Z$$
l'homomorphisme qui envoie la famille $\{y_w\}_{w \in \Omega_L}$
sur $\sum_{w \in \Omega_L} [\F_w:\F]w(y_w)$.
Cet homomorphisme est trivial sur l'image diagonale de $L^{\times}$
dans $\I_L$ (loi de r\'eciprocit\'e, ``le nombre des z\'eros est
\'egal au nombre des p\^oles"). Il est aussi trivial sur
le sous-groupe compact maximal $\prod_{w\in \Omega_L} O_w^{\times}$.

Si l'extension $L/K$ est de plus galoisienne de groupe $G$,
le groupe $G$ agit naturellement sur $\I_L$, et trivialement
sur $\Z$.
On v\'erifie que l'homomorphisme $\deg_{L,\F} : \I_L \to \Z$
est $G$-\'equivariant.

Soient $\F,K,L$ comme ci-dessus, avec l'extension $L/K$ galoisienne
de groupe $G$. Soit $T$ un $K$-tore d\'eploy\'e par $L$.
L'homomorphisme $\deg_{L,\F} : \I_L \to \Z$ induit un $G$-homomorphisme
 $$\deg_{L,\F,T} : T_* \otimes \I_L \to T_*,$$
qui est $G$-\'equivariant, l'action de $G$ \`a gauche \'etant
l'action simultan\'ee sur $T_*$ et $\I_L$.
L'homomorphisme  ainsi obtenu est fonctoriel en les $K$-tores 
d\'eploy\'es par $L$. Il est nul sur $T_*\otimes L^{\times}$
et sur $T_* \otimes (\prod_{w\in \Omega_L} O_w^{\times})$.

Cet homomorphisme induit sur les points fixes sous $G$
un homomorphisme de groupes ab\'eliens
 $$\deg_{L,\F,T} : T(\A_K)=(T_* \otimes \I_L)^G \to T_*^G,$$
c'est-\`a-dire des id\`eles de $T$ (sur $K$) vers $T_*^G$.
Cet homomorphisme s'annule  sur le
sous-groupe compact maximal des id\`eles de $T$, qui est $\prod_{v \in
\Omega_K} T(O_v)=(T_*
\otimes (\prod_{w\in
\Omega_L} O_w^{\times}))^G$ et sur $T(K) \subset T(\A_K)$.
  L'homomorphisme  ainsi obtenu est fonctoriel
en les $K$-tores  d\'eploy\'es par le corps $L$.

Soit 
$$0 \to F_* \to P_* \to T_* \to 0 \hskip3cm $$
une suite exacte de $G$-r\'eseaux du type (1) (r\'esolution
coflasque de $T_*$). 

\medskip

{\bf Question globale } {\it L'application compos\'ee de 
$\deg_{L,\F,P} : P(\A_K) = (P_* \otimes \I_L)^G \to P_*^G$
et de $P_*^G \to T_*^G$
a-t-elle m\^eme image
que l'application $\deg_{L,\F,T} : T(\A_K)= (T_* \otimes \I_L)^G
\to T_*^G$ ?}

\medskip
(Par fonctorialit\'e, la premi\`ere image est contenue dans la
seconde.)

\medskip

On voit imm\'ediatement que  la r\'eponse \`a cette question ne d\'epend pas du choix du corps fini $\F \subset K$. En utilisant les propri\'et\'es des r\'esolutions flasques et coflasques, on voit aussi que
la r\'eponse \`a cette question ne d\'epend que du $K$-tore $T$, elle ne d\'epend ni
du choix du corps de d\'eploiement $L/K$ ni du choix de la
 r\'esolution coflasque (1) de $T_{*}$.

\bigskip

Montrons  que cette question est \'equivalente \`a celle rencontr\'ee par Bourqui dans [B].
Soit $G$ un groupe fini et $M$ un $G$-r\'eseau.
On note $M^{ \circ} $ le $G$-r\'eseau $\hom(M,\Z)=\hom_{\bf Z}(M,\Z)$.

\medskip

{\bf Lemme 3.1} {\it  Soit $M$ un $G$-r\'eseau. L'inclusion
$M^G \subset M$ induit une application injective
$(M^{\circ })^G \hookrightarrow (M^G)^{\circ }$
\`a conoyau fini.}

{\it D\'emonstration}
Consid\'erons la suite exacte
$$0 \to M^G \to M \to R \to 0$$
d\'efinissant $R$ comme le conoyau de l'inclusion
naturelle. Le groupe ab\'elien $R$ est sans torsion,
la suite est donc scind\'ee comme suite de groupes ab\'eliens.
On a donc la suite exacte de $G$-modules duale
$$ 0 \to R^{\circ } \to M^{\circ } \to (M^G)^{\circ } \to 0.$$
Soit $\varphi \in (R^{\circ})^G={\rm Hom}_G(R,\Z)$.
Pour tout $r \in R$, on a $N_{G}r=0$.
On a donc $0=\varphi(N_{G} r)=N_{G}(\varphi (r))= n (\varphi (r))$, o\`u $n>0$ est l'ordre de $G$.
Ainsi $\varphi (r)=0$ pour tout $r \in R$, i.e. $\varphi=0$. 
Ceci \'etablit  $(R^{\circ})^G=0$.
Le d\'ebut de la suite exacte de $G$-cohomologie associ\'ee \`a la derni\`ere
suite exacte s'\'ecrit donc
$$ 0 \to  (M^{\circ })^G \to (M^G)^{\circ } \to H^1(G,R^{\circ }),$$
ce qui \'etablit le lemme. \cqfd

\bigskip

Soit $\chi \in T^*{^G}$, c'est-\`a-dire un caract\`ere, d\'efini
sur $K$, du $K$-tore  $T$. La donn\'ee d'un
tel \'el\'ement $\chi$ \'equivaut \`a celle d'un homomorphisme
$G$-\'equivariant $T_* \to \Z$. Celui-ci induit un
homomorphisme
$G$-\'equivariant
$T_* \otimes \I_L \to \I_L$ et donc, en prenant les points fixes sous $G$,
un homomorphisme
$T(\A_K) \to \I_K$. On peut composer ceci avec l'application 
$\deg_{K,\F} : \I_K \to \Z$. On d\'efinit ainsi une application
bilin\'eaire 
$$ T(\A_K) \times T^*{^G} \to \Z$$
soit encore 
$$ \deg_{T,K,\F} : T(\A_K) \to \hom( T^*{^G}, \Z),$$
qui est nulle sur l'image de $T(K)$ dans $T(\A_K)$ et
 sur tout
\'el\'ement de $\prod_{v \in \Omega}T(O_v)$.
L'application ainsi d\'efinie est fonctorielle en le $K$-tore $T$.
Elle ne d\'epend pas du choix du corps de d\'eploiement $L/K$ de $T$.
Lorsque
le corps $\F$ est alg\'ebriquement ferm\'e dans $K$, elle
co\"{\i}ncide avec l'application  ${\rm deg}_{T}$ d\'efinie  au \S 2.3 de [B].
Un calcul analogue au calcul local fait au \S 2 (formule $\phi_{K,L}=e\psi_{K}$)
\'etablit le lemme suivant.
\medskip

{\bf Lemme 3.2} {\it La fl\`eche compos\'ee 
$$T(\A_K) \to T_*^G \hookrightarrow  \hom( T^*{^G}, \Z),$$
o\`u la premi\`ere fl\`eche est $\deg_{L,\F,T} $
et la seconde fl\`eche  l'inclusion naturelle du Lemme 3.1,
est \'egale \`a $[L:K]. \deg_{T,K,\F}$. \cqfd}

\bigskip

En appliquant le foncteur $\hom( \bullet, \Z)$
\`a la suite (1), on obtient 
 une suite exacte de $G$-r\'eseaux
$$ 0 \to T^* \to P^* \to F^* \to 0,$$
une fl\`eche $T^*{^G} \to P^*{^G} $ 
et une fl\`eche 
$$\hom( P^*{^G}, \Z) \to \hom( T^*{^G}, \Z).$$

 L'application compos\'ee
$$ \deg_{P,K,\F} : P(\A_K) \to 
 \hom( P^*{^G}, \Z) \to \hom( T^*{^G}, \Z)
$$ 
a son image contenue dans 
celle de l'application 
$$ \deg_{T,K,\F} : T(\A_K) \to 
  \hom( T^*{^G}, \Z).$$ 

Le probl\`eme rencontr\'e dans [B] est le suivant : {\it Ces deux images
co\"{\i}ncident-elles ?} Lorsque $\F$ est alg\'ebriquement ferm\'e dans $K$,
le quotient des deux images  est le groupe fini not\'e ${\cal K}_{T}$ dans [B] (\S 2.7).
Le Lemme 3.2 montre que le probl\`eme
se traduit imm\'ediatement en la question globale.

\bigskip

Comme le note Bourqui ([B], Prop. 2.15), il est des cas o\`u la question globale a une r\'eponse
positive. 

(a) C'est le cas si  le $K$-tore $T$ est anisotrope, car alors
 $T_*^G=0$ et $T^*{^G}=0$. 

(b) C'est le cas si  le  corps des constantes de $K$ 
co\"{\i}ncide avec celui de $L$,  i.e.  est alg\'ebriquement ferm\'e dans  $L$.
Bourqui montre ([B], \S 2.9, Lemme 2.18)
que sous l'hypoth\`ese que le corps $\F$ est alg\'ebriquement ferm\'e
dans $K$ et dans $L$, l'application
$\deg_{T,K,\F} : T(\A_K) \to 
  \hom( T^*{^G}, \Z)$ envoie 
$N_G T(\A_L) \subset T(\A_K)$ surjectivement sur $\hom( T^*{^G}, \Z)$.
Comme par ailleurs $P(\A_L)$ se surjecte sur $T(\A_L)$,
ceci \'etablit la surjectivit\'e voulue.

(c) C'est le cas si le $K$-tore $T$ satisfait l'approximation
faible   ([B], Lemme 2.13 et Prop. 2.15).
A ce sujet, on a l'\'enonc\'e plus g\'en\'eral suivant.

\bigskip

{\bf Proposition 3.3} {\it  Soit $T$ un $K$-tore d\'eploy\'e par l'extension galoisienne finie
$L/K$ de groupe de Galois $G$. Soit $S$ l'ensemble fini des places de $K$
telles que le groupe de Galois local ne soit pas m\'etacyclique. Si la r\'eponse
\`a la question semi-locale pour $T$ et $S$ est affirmative, i.e. si le sous-groupe 
 ouvert $T(K).\prod_{v \in S} T(O_{v})$ de $\prod_{v \in S} T(K_{v})$
co\"{\i}ncide  avec $\prod_{v \in S} T(K_{v})$, alors la question globale pour $T$
a une r\'eponse affirmative. }

\medskip

{\it D\'emonstration} Soit 
$$1 \to F \to P \to T \to 1$$
une r\'esolution flasque de $T$ par des $K$-tores d\'eploy\'es par $L$.
En toute place $v$ de $K$ non dans $S$,  le th\'eor\`eme d'Endo et Miyata rappel\'e au \S 1
assure que le $K_{v}$-tore $F\times_{K}K_{v}$ est un facteur direct d'un $K_{v}$-tore
quasitrivial, ce qui implique $H^1(K_{v},F)=0$, et donc $P(K_{v}) \to T(K_{v})$ surjectif.
Soit $\xi=\{t_v\}_{v \in \Omega} \in T(\A_K)$. 
Si la question semi-locale pour $T$ et $S$ a une r\'eponse affirmative, il existe
$t \in T(K)$ tel que toute composante de $ t.\xi$ pour $v \in S$ soit
dans $T(O_{v})$. En toute place $v$ non dans $S$, 
la composante de $ t.\xi$ est dans l'image de $P(K_{v})$.
On voit ainsi que $t.\xi$ est le produit d'un id\`ele appartenant
\`a $\prod_{v\in \Omega} T(O_{v})$ et d'un id\`ele dans l'image
de $P(\A_{K}) \to T(\A_{K})$. L'application 
$ \deg_{T,K,\F} : T(\A_K) \to 
  \hom( T^*{^G}, \Z)$ est nulle sur le groupe compact $\prod_{v\in \Omega} T(O_{v})$,
  et par r\'eciprocit\'e elle est nulle sur $T(K) \subset T(\A_{K})$. 
  On voit donc que l'image de $\xi$ dans $\hom( T^*{^G}, \Z)$
  est l'image d'un \'el\'ement de l'application compos\'ee
  $P(\A_{K}) \to T(\A_{K}) \to \hom( T^*{^G}, \Z).$  \cqfd

\medskip

(d) Le rapporteur note que l'on peut aussi, dans le cas global ici consid\'er\'e, \'etablir l'analogue du cas (iv) de la Proposition 2.1.

\bigskip

{\bf \S 4. Construction et \'etude d'un r\'eseau muni d'une action du
groupe de Klein.}

\medskip

 Etant donn\'e un groupe fini $G$, on note $\Z$ le 
r\'eseau $\Z$ avec action triviale de $G$ et on note
$\Z[G]$ le $G$-r\'eseau
standard de $\Z$-base les \'el\'ements de $G$.
On note $I_G$ le noyau de l'augmentation $\varepsilon_G : \Z[G] \to \Z$.
On note $N_G= \sum_{g \in G}g \in \Z[G]$.
On sait ([CTS1], Prop. 15 p. 206)
que si les \'el\'ements $\sigma_i, i \in I,$ 
engendrent le groupe $G$, alors le $G$-homomorphisme
$\oplus_{i \in I} \Z[G] \to I_G$
qui sur la coordonn\'ee $i$ envoie 1 sur $1-\sigma_i$
est surjectif et a pour noyau un $G$-module coflasque.

Soit $G = < \sigma,\tau>$ avec $\sigma^2=1, \tau^2=1, \sigma \tau=
\tau
\sigma$. Soit  $T_*$ le $G$-r\'eseau noyau de
l'homo\-mor\-phisme
$$ \Z[G] \oplus I_G  \to \Z[G] \oplus \Z[G]$$
donn\'e par
$$(t,x) \mapsto (\sigma t-t-x-\tau x,\tau t -t - x -\sigma x).$$

\medskip

{\bf Lemme 4.1} {\it 

(i) L'application compos\'ee de l'inclusion $T_{*} \subset \Z[G] \oplus I_G  $
et de la projection $\Z[G] \oplus I_G   \to \Z[G]$ est injective.

(ii) On a la suite exacte de $G$-r\'eseaux
$$0 \to \Z \to T_* \to I_G \to 0,$$
o\`u l'application $\Z \to T_* \hookrightarrow \Z[G] \oplus I_G$
est donn\'ee par $1 \mapsto (N_G,0)$
et o\`u l'application $T_* \to I_G$ est
la compos\'ee de $T_* \hookrightarrow \Z[G] \oplus I_G$
et de la projection sur le second facteur.}

\medskip

La preuve est laiss\'ee au lecteur. Bien que nous n'en ayons
pas besoin, notons qu'on a  une suite exacte longue
$$ 0 \to T_* \to \Z[G] \oplus I_G \to (1-\sigma)(\Z+\Z\tau) \oplus
(1-\tau)(\Z+\Z\sigma) \to \Z \to 0.$$
Dans cette suite, chacun des $G$-modules 
 $(1-\sigma)(\Z+\Z\tau)$
et
$(1-\tau)(\Z+\Z\sigma)$ est un $G$-sous-module de $\Z[G]$,
la fl\`eche compos\'ee
$$\Z[G] \oplus I_G \to (1-\sigma)(\Z+\Z\tau) \oplus
(1-\tau)(\Z+\Z\sigma) \subset \Z[G] \oplus \Z[G]$$
est la fl\`eche 
$(t,x) \mapsto (\sigma t-t-x-\tau x,\tau t -t - x -\sigma x)$
dont le noyau d\'efinit $T_*$, l'application 
$$(1-\sigma)(\Z+\Z\tau) \oplus
(1-\tau)(\Z+\Z\sigma)  \to \Z$$
envoie
$((1-\sigma)(a+b\tau),
(1-\tau)(c+d\sigma))$ sur $a+c-b-d$. 

\medskip

L'homomorphisme $\varphi : \Z[G] \oplus \Z[G] \to I_G$ 
donn\'e 
par $ (a,b) \mapsto a(1-\sigma)+b(1-\tau)$
donne une r\'esolution coflasque 
$$0 \to F_* \to \Z[G] \oplus \Z[G] \to I_G \to 0$$
de $I_G$.
L'image r\'eciproque 
$$0 \to \Z \to P_* \to \Z[G] \oplus \Z[G] \to 0$$
de la suite du Lemme 4.1
par $\varphi$ d\'efinit une extension
de $\Z[G] \oplus \Z[G] $ par $\Z$. Toute telle
extension de $G$-r\'eseaux est scind\'ee ($H^1(G,\Z)=0$).
Il existe donc un $G$-rel\`evement
$\Z[G] \oplus \Z[G] \to T_*$ de $\varphi$.
On peut prendre 
pour ce rel\`evement la fl\`eche
qui \`a $(1,0)$ associe $(\sigma+\sigma\tau, 1-\sigma)$
et \`a $(0,1)$ associe $(\tau + \sigma \tau, 1-\tau)$.
(Deux tels rel\`evements diff\`erent par
 une application $(a,b) \mapsto (an+bm)(N_G,0)$,
o\`u $n$ et $m$ peuvent \^etre pris arbitraires.)

\medskip

L'image r\'eciproque de la r\'esolution coflasque
$$0 \to F_* \to \Z[G] \oplus \Z[G] \to I_G \to 0$$
par la fl\`eche $T_* \to I_G$ est une
r\'esolution coflasque de $T_*$:
$$ 0 \to F_* \to P_* \to T_* \to 0,$$
o\`u $P_*=\Z \oplus \Z[G] \oplus \Z[G],$
la fl\`eche compos\'ee
$$\Z \oplus \Z[G] \oplus \Z[G] = P_* \to T_* \hookrightarrow \Z[G] \oplus
I_G \to \Z[G] $$
(la derni\`ere fl\`eche \'etant la projection sur le facteur 
$\Z[G]$ de $\Z[G] \oplus I_G$)
\'etant donn\'ee par
$$(a,b,c) \mapsto N_G a +(\sigma+\sigma \tau)b + 
(\tau + \sigma\tau)c.$$

\bigskip

{\bf \S 5. La question locale a une r\'eponse n\'egative.}

\medskip

Soit $K$ un corps local de corps r\'esiduel le corps fini
$\F$, suppos\'e de caract\'eristique $p \neq 2$.  Soit $v : K^{\times}
\to \Z$ la valuation normalis\'ee.
Le groupe $H^1(K,\Z/2)$ est alors isomorphe \`a $(\Z/2)^2$.
Soit $L/K$ l'unique extension galoisienne  de $K$ de groupe
$G \simeq \Z/2 \times \Z/2$. Soit $w$ la valuation normalis\'ee
sur
$L$. L'extension $L/K$ est ramifi\'ee, l'indice de ramification est 2,
pour
$\alpha \in K^{\times}$, on a $w(\alpha)=2v(\alpha)$.
Le corps r\'esiduel $\F_w$ de $L$ est une extension
quadratique de $\F$. 

Soit $\pi $ une uniformisante de $K$, et soit $u \in K^{\times}$
une unit\'e qui n'est pas un carr\'e. Soient $L_1=K({\sqrt \pi }) \subset
L$ et $L_2=K({\sqrt u}) \subset L$. Les sous-extensions $L_1/K$
et $L_2/K$ de $L/K$ sont quadratiques. Appelons  $\sigma \in G$
l'\'el\'ement non trivial
fixant $L_1$ et $\tau \in G$ l'\'el\'ement non trivial fixant $L_2$.

La r\'esolution coflasque de $G$-modules
$$ 0 \to F_* \to P_* \to T_* \to 0$$
consid\'er\'ee au \S 4
induit sur les $K$-points des $K$-tores 
d\'eploy\'es par $L$ associ\'es
un homomorphisme
$K^{\times} \times L^{\times} \times L^{\times} = P(K) \to T(K) \subset
L^{\times}
\times L^{{\times},1},$ o\`u $L^{{\times},1} \subset L^{\times}$ est le
sous-groupe des
\'el\'ements de norme 1. L'image de $P(K)$ dans $T(K)$ est le sous-groupe
$RT(K)$ des \'el\'ements $R$-\'equivalents \`a 1.
La compos\'ee de l'inclusion  $T(K) \subset
L^{\times}
\times L^{{\times},1}$ et de la projection sur le premier facteur de ce dernier produit
d\'efinit un plongement $T(K) \subset L^{\times}$.

L'homomorphisme compos\'e
$K^{\times} \times L^{\times} \times L^{\times} = P(K) \to T(K) \subset
L^{\times}$ est alors donn\'e par $$(\alpha,\beta,\gamma) \mapsto
\alpha.((1+\tau)\sigma
\beta). ((1+\sigma) \tau \gamma) \in L^{\times}$$
La valuation normalis\'ee $w$ d'un tel \'el\'ement de $L^{\times}$ est
paire. En effet pour tout \'el\'ement  $\alpha \in K^{\times} \subset
L^{\times}$ on a
$w(\alpha)=2v(\alpha)$ et pour tout \'el\'ement $\delta \in L^{\times}$
et tout
$g \in G$, on a $w(\delta.g(\delta))=2w(\delta)$.

Supposons que $-1$ est un carr\'e dans $\F$, donc dans $K$.
Soit $i \in K$ tel que $i^2=-1$.

Consid\'erons l'\'el\'ement $({\sqrt{u\pi }}, i) \in L^{\times} \times
L^{\times}$. Comme $i \in K$, on a $N_G(i)=1$, donc 
$({\sqrt{u\pi }}, i) \in L^{\times} \times L^{{\times},1}$.
Par ailleurs $\sigma({\sqrt{u\pi }})/{\sqrt{u\pi }}=-1= i \tau(i)$
et $\tau( {\sqrt{u\pi }} )/{\sqrt{u\pi }}=-1= i \sigma(i)$.
Donc $({\sqrt{u\pi }}, i)$ appartient \`a $T(K) \subset L^{\times} \times
L^{{\times},1}$. L'image de $({\sqrt{u\pi }}, i)$ par la projection sur le
premier facteur est ${\sqrt{u\pi }} \in L^{\times}$, dont la valuation
normalis\'ee est 1.

Comme l'image du sous-groupe compact maximal de $T(K)$
dans $L^{\times}$ est dans $O_L^{\times}$, on conclut que, via la
projection
$T(K) \to L^{\times}$ donn\'ee par le premier facteur, l'image du
sous-groupe engendr\'e par le sous-groupe compact maximal de $T(K)$ et
le sous-groupe image de $P(K)$ consiste en des \'el\'ements de
valuation normalis\'ee paire de $L^{\times}$, alors qu'il existe
un \'el\'ement de $T(K)$ dont l'image dans $L^{\times}$ est de valuation
normalis\'ee 1. Ceci \'etablit $T(K)\neq T(O_K). RT(K)$. \cqfd

\medskip

D'apr\`es la Proposition 2.2, la r\'eponse n\'egative \`a la question semi-locale,
que nous allons donner au \S 6, donne aussi une r\'eponse n\'egative
\`a la question  locale. Mais d'une part le calcul du pr\'esent paragraphe
est utilis\'e au \S 6, d'autre part il vaut
aussi pour un corps local $K$ de caract\'eristique nulle.

\bigskip
 
{\bf \S 6. Les questions semi-locale et globale  ont  une r\'eponse n\'egative.}

\medskip

Soit $\F$ un corps fini de caract\'eristique impaire,
tel que $-1$ soit un carr\'e dans $\F$. Soit $\F'$
l'extension quadratique de $\F$.
Soient $K=\F(\lambda)$ le corps des fractions rationnelles sur 
$\F$, puis $L_1=\F({\sqrt \lambda})$, $L_2=\F'(\lambda)$  et
$L=\F'({\sqrt \lambda})$.

Le groupe de Galois $G$ de $L/K$ est $\Z/2 \times \Z/2
=<\sigma,\tau>$, o\`u $\sigma$ est l'\'el\'ement non trivial qui laisse
fixe $L_1$ et $\tau$ l'\'el\'ement non trivial qui laisse fixe $L_2$.

Soit $T_*$   comme au \S 4. On dispose donc
du $G$-homomorphisme $$\Z \oplus \Z[G] \oplus \Z[G] = P_* \to T_*$$
qui, compos\'e avec la projection de $T_* \subset \Z[G] \oplus I_G$
sur le premier facteur, se lit
$$(a,b,c) \mapsto N_G a +(\sigma+\sigma \tau)b + 
(\tau + \sigma\tau)c \in \Z[G].$$

Comme au \S 3, on note $\deg_{L,\F} : \I_L \to \Z$
le degr\'e relatif au corps de base $\F$.

\medskip

{\bf Proposition 6.1} {\it L'image de l'application compos\'ee de 
$\deg_{L,\F,P} : P(\A_K) = (P_* \otimes \I_L)^G \to P_*^{G}$
et de $P_*^{G} \to T_*^G$
est strictement contenue dans l'image
 de l'application $\deg_{L,\F,T} : T(\A_K)= (T_* \otimes \I_L)^G
\to T_*^G$.}

\medskip
{\it D\'emonstration} Pour \'etablir cette proposition, il suffit de
montrer les deux faits suivants :

(a)  Consid\'erons l'application compos\'ee
$$P(\A_K) = (P_* \otimes \I_L)^G \to P_*^{G} \to T_*^G \to \Z[G]^G \simeq
\Z,$$
o\`u le dernier isomorphisme $\Z[G] \simeq \Z$  est l'inverse de 
l'application envoyant $1$ sur $N_G$, et o\`u la fl\`eche
$ T_*^G \to \Z[G]^G$ est induite par la projection de $T_*
\subset \Z[G] \oplus I_G$ sur le premier facteur.
L'image de cette application compos\'ee
 est contenue dans $4\Z$.

(b) L'image de l'application compos\'ee 
$$T(\A_K)=(T_* \otimes \I_L)^G \to T_*^G \to \Z[G]^G \simeq
\Z,$$
contient  $2 \in \Z$.

\medskip
Etablissons le point (a). L'application consid\'er\'ee  est
obtenue de la fa\c con suivante. On consid\`ere le $G$-homomorphisme
$\Z \oplus \Z[G] \oplus \Z[G]  \to \Z[G]$
donn\'e par $$(a,b,c) \mapsto N_G a +(1+\tau)\sigma b + 
(1 + \sigma)\tau c$$ et l'homomorphisme
 $\deg_{L,\F} : \I_L \to \Z$. On tensorise
$$(\Z \oplus \Z[G] \oplus \Z[G]) \otimes \I_L \to \Z[G]\otimes \I_L \to
\Z[G]
$$
 et on prend les points fixes sous $G$. Ceci donne
$$ \I_K \oplus (\Z[G] \otimes \I_L)^G \oplus (\Z[G] \otimes \I_L)^G
\to (\Z[G]\otimes \I_L)^G \to \Z[G]^G = \Z N_G,$$
qu'on identifie 
\`a
$$ \I_K \oplus \I_L \oplus \I_L \to \I_L \to \Z,$$
o\`u $\I_K \to \I_L$ est l'inclusion naturelle,
la premi\`ere application $\I_L \to \I_L$ est donn\'ee
par $\xi \mapsto (1+\tau)\sigma\xi$, la seconde par
$\eta\mapsto (1 + \sigma)\tau\eta$, et la fl\`eche
$\I_L \to \Z$ est $\deg_{L,\F}$.
La compos\'ee de l'application diagonale $\I_K \to \I_L$
et de $\deg_{L,\F} : \I_L \to \Z$ est $4.\deg_{K,\F}$. 
L'application degr\'e $\deg_{L,\F} :  \I_L \to \Z$ est
$G$-\'equivariante. Ainsi l'image de $(1+\tau)\sigma\xi \in \I_L$
dans $\Z$ appartient \`a $2 \deg_{L,\F}(\I_L) \subset \Z$.
Mais pour tout compl\'et\'e $L_w$ de $L$ le corps r\'esiduel $\F_w$
contient $\F'$, et donc $[\F_w:\F]$ est pair. 
On a donc $\deg_{L,\F}(\I_L) \subset 2.\Z$, et 
l'image de $(1+\tau)\sigma\xi \in \I_L$
dans $\Z$ est dans $4.\Z$. L'argument est le m\^eme pour
l'image de $(1 + \sigma)\tau\eta$.

\medskip

Pour \'etablir le point b), consid\'erons simplement
la place $v$ de $K$ d\'efinie par $t=0$. Il y a une seule
place $w$ de $L$ au-dessus de $v$, et l'extension locale
$L_w/K_v$ est $\F'(({\sqrt \lambda}))/\F((\lambda))$, elle est du type
consid\'er\'e au \S 5. L'application compos\'ee
$$T(K_v) \subset T(\A_K) \to T_*^G \to \Z[G]^G \simeq
\Z$$
envoie $(t,x) \in T(K_v) \subset L_w^{\times} \times L_w^{{\times},1}$
sur $[\F_w:\F_v] w(t) = [\F':\F] w(t)=2 w(t)$, o\`u $w$ est la valuation
normalis\'ee de $L_w$. On a vu au \S 5
qu'il existe un \'el\'ement $(t,x) \in T(K_v)$
avec $w(t)=1$. Ceci ach\`eve la d\'emonstration.\cqfd

Ainsi la question globale (\S 3) a une r\'eponse n\'egative. D'apr\`es la
Proposition 3.3, ceci implique que la question semi-locale (sur un corps
de fonctions d'une variable sur un corps fini) a une r\'eponse
n\'egative.

\bigskip

{\it Remerciements}

Ce travail 
a \'et\'e entrepris en mai 2004, lors d'un s\'ejour du second
auteur (V.S.) au laboratoire de Math\'ematiques de l'Universit\'e
Paris-Sud. Ce s\'ejour a \'et\'e rendu possible
gr\^ace au soutien du Centre
 franco-indien pour la recherche avanc\'ee (CEFIPRA, IFCPAR),
projet num\'ero 2501-1. 

\bigskip

{\bf Bibliographie}

\bigskip

[BT] Victor V. Batyrev and Yuri Tschinkel,
Rational points of bounded height on compactifications of anisotropic tori,
Internat. Math. Research Notices {\bf 12} (1995) 591--635.

[B] D. Bourqui, Constante de Peyre des vari\'et\'es
toriques en caract\'eristique positive, pr\'epublication 2004, disponible sur le serveur
arXiv sous la r\'ef\'erence math.NT/0501409.

[CTS1]  J.-L. Colliot-Th\'el\`ene et J.-J. Sansuc, La R-\'equivalence
sur les tores, Ann. scient.  \'Ec. Norm. Sup. {\bf 10} (1977) 175--229.

[CTS2] J.-L. Colliot-Th\'el\`ene et J.-J. Sansuc,
Principal homogeneous spaces under flasque tori : Applications, J. Algebra
{\bf 106} (1987) 148--205.

[P1] E. Peyre,  Hauteurs et mesures de Tamagawa sur les vari\'et\'es de Fano, Duke Math. J. {\bf 79} (1995)  101--218.

[P2] E. Peyre, Points de hauteur born\'ee, topologie ad\'elique et mesure de Tamagawa,
Journal de  Th\'eorie des Nombres de Bordeaux {\bf 15} (2003) 319--348.

\bigskip

\vskip1cm

Jean-Louis Colliot-Th\'el\`ene

C.N.R.S., Math\'ematiques,

UMR 8628,

B\^atiment 425,

Universit\'e Paris-Sud,

F-91405 Orsay

FRANCE

colliot@math.u-psud.fr

\vskip1cm

Venapally Suresh,

Department of Mathematics and Statistics,

University of Hyderabad, 

P.O. Central University, 

Gachibowli, 

Hyderabad 500 046, 

Andhra Pradesh, 

INDE

 vssm@uohyd.ernet.in

\bye